\documentclass{article}
\usepackage{amssymb}
\usepackage{graphicx}
\usepackage{amsmath}

\input{tcilatex}
\begin{document}

\begin{center}
{\LARGE Higher secants of spinor varieties\bigskip \medskip }

{\Large Elena Angelini\medskip }

{\scriptsize Universit\`{a} di Firenze, Dipartimento di Matematica "Ulisse
Dini",}

{\scriptsize Viale G. B. Morgagni 67/a, 50134 Firenze, Italy}{\Large %
\bigskip }
\end{center}

\begin{quotation}
\noindent {\small ABSTRACT. Let }$S_{h}$ be the even pure spinors variety of
a complex vector space $V$ of even dimension $2h$ endowed with a non
degenerate quadratic form $Q$ and let $\sigma _{k}\left( S_{h}\right) $ be
the $k$-secant variety of $S_{h}$. We decribe a probabilistic algorithm
which computes the complex dimension of $\sigma _{k}\left( S_{h}\right) $.
Then, by using an inductive argument, we get our main result: $\sigma
_{3}\left( S_{h}\right) $ has the expected dimension except when $h\in
\left\{ 7,8\right\} $. Also we provide theoretical arguments which prove
that $S_{7}$ has a defective $3$-secant variety and $S_{8}$ has defective $3$%
-secant and $4$-secant varieties.\footnote{%
}$\bigskip $
\end{quotation}

\section{Introduction\emph{\protect\medskip }}

\noindent In this paper we study the higher secant varieties of spinor
varieties.

\noindent We consider a complex $2h$-dimensional vector space $V$ and a non
degenerate quadratic form $Q$ defined on it. The \textit{space of spinors }%
associated to $\left( V,Q\right) $ can be identified with the space of the
spin representation of $Cl\left( V,Q\right) $, the \textit{Clifford algebra}
generated by $V$. In particular, \textit{pure spinors }represent, from a
geometrical point of view, the set of all maximal totally isotropic vector
subspaces of $V$, which is a projective variety, called \textit{spinor
variety}. For simplicity, we consider one of its two irreducible isomorphic
components, i.e. the \textit{even pure spinors variety}, which we denote by $%
S_{h}$.

\noindent Let $X$ be a non-degenerate projective variety in $\mathbb{P}%
^{N}\left( \mathbb{C}\right) $; then $\sigma _{k}\left( X\right) $ indicates
the $k$\textit{-secant variety }of $X$, that is the Zariski closure of the
union of all linear spaces spanned by $k$ points of $X$, see $\left( \cite{Z}%
\right) $ and $\left( \cite{L}\right) $ for several applications. It's easy
to check the following inequality:%
\begin{equation*}
\dim _{\mathbb{C}}\sigma _{k}\left( X\right) \leq \min \left\{ k\dim _{%
\mathbb{C}}X+k-1,N\right\} \text{.}
\end{equation*}%
If the equality holds, then we say that $\sigma _{k}\left( X\right) $ has
the \textit{expected dimension}, otherwise $X$ is said to be $k$-\textit{%
defective} and%
\begin{equation*}
\delta _{k}=\min \left\{ k\dim _{\mathbb{C}}X+k-1,N\right\} -\dim _{\mathbb{C%
}}\sigma _{k}\left( X\right)
\end{equation*}%
is its $k$-\textit{defect}. The problem of determining the complex dimension
of $\sigma _{k}\left( X\right) $ is called the \textit{defectivity problem }%
for $X$. If $\nu _{d}\left( \mathbb{P}^{n}\left( \mathbb{C}\right) \right) $
is the Veronese variety then $\sigma _{k}\left( \nu _{d}\left( \mathbb{P}%
^{n}\left( \mathbb{C}\right) \right) \right) $ has the expected dimension
except in some particular cases, $\left( \cite{A-H}\right) $, $\left( \cite%
{B-O}\right) $. Concerning Grassmannians and Segre varieties, this problem
has been studied by several authors but it's still open, as we can see,
respectively, in $\cite{B-D-DG}$ and $\cite{A1}$; for related results see
also $\cite{A}$, $\cite{B-B}$ and $\cite{C-G-G}$. At the best of my
knowledge, the case of spinor varieties is almost absent in the mathematical
literature; it's known that $\sigma _{2}\left( S_{h}\right) $ has always the
expected dimension $\left( \cite{K}\right) $, but for $k\geq 3$ the problem
was completely open.

\noindent By using \textit{Macaulay2} software system, we construct a
probabilistic algorithm which allow us to compute the dimension of $\sigma
_{k}\left( S_{h}\right) $ by studying the span of the tangent spaces at $k$
chosen random points, for $h\leq 12$. Afterwards, by using induction, we get
our main result:\medskip

\noindent \textbf{Theorem 1.1} (\textit{i}) $\sigma _{3}\left( S_{h}\right) $
\textit{has the expected dimension}, \textit{except when} $h\in \left\{
7,8\right\} $.\smallskip

\noindent (\textit{ii}) $S_{7}$ \textit{has a defective} $3$\textit{-secant
variety and} $S_{8}$ \textit{has defective} $3$\textit{-secant and} $4$%
\textit{-secant varieties. In particular }$\dim _{\mathbb{C}}\sigma
_{3}\left( S_{7}\right) =58$, $\dim _{\mathbb{C}}\sigma _{3}\left(
S_{8}\right) =85$ \textit{and} $\dim _{\mathbb{C}}\sigma _{4}\left(
S_{8}\right) =111$.\medskip

\noindent We remark that the main tool of our investigation is the
parametrization of $S_{h}$ with all \textit{principal sub-Pfaffians} of a
skew symmetric matrix of size $h$.

\noindent The paper is organized in six sections. In the second one we
introduce Clifford algebras and spinor varieties, following $\cite{C}$, $%
\cite{M}$ and $\cite{Ang}$; in the third we recall the main definitions and
properties of higher secant varieties, $\left( \cite{L}\right) $, $\left( 
\cite{Z}\right) $. Finally, sections four, five and six are devoted to our
main results.\smallskip

\noindent This article is based upon the author's laurea thesis and the main
result confirms its final conjectures, $\left( \cite{Ang}\right) $. Thanks
are due especially to Giorgio Ottaviani for his guidance and insight.

\section{Clifford algebras and spinors$\protect\medskip $}

\noindent Let $V$ be a vector space over $\mathbb{C}$ of even dimension $%
n=2h>0$. Let $Q$ be a quadratic form on $V$ such that the corresponding
symmetric bilinear form $B$ is non degenerate.\smallskip

\noindent We denote by $Cl\left( V,Q\right) =T\left( V\right) /I_{Q}\left(
V\right) $ the \textit{Clifford algebra associated to} $\left( V,Q\right) $,
where $T\left( V\right) $ is the tensor algebra of $V$ and $I_{Q}\left(
V\right) \subset T\left( V\right) $ is the two-sided ideal generated by the
elements%
\begin{equation*}
v\otimes v-Q\left( v\right) \cdot 1
\end{equation*}%
with $v\in V$.\smallskip

\noindent Let%
\begin{equation*}
Cl(V,Q)_{\pm }=T\left( V\right) _{\pm }/I_{Q}\left( V\right) \cap T\left(
V\right) _{\pm }
\end{equation*}%
where\textit{\ }$T\left( V\right) _{+}$ and $T\left( V\right) _{-}$ denote
the set of even and odd tensors, respectively. Then $Cl(V,Q)_{+}$ is a
subalgebra of $Cl\left( V,Q\right) $ and%
\begin{equation*}
Cl\left( V,Q\right) =Cl(V,Q)_{+}\oplus Cl(V,Q)_{-}\text{.}
\end{equation*}%
In particular, we call \textit{even}\ the elements of $Cl(V,Q)_{+}$ and 
\textit{odd}\ those of $Cl(V,Q)_{-}$.\medskip

\noindent Let $E$ and $F$ be maximal totally isotropic vector subspaces of $V
$ such that $V=E\oplus F$, let $f$ be the product in $Cl(V,Q)$ of the
elements of a basis of $F$. It can be proved $\left( \cite{C}\right) $ that
there's only one irreducible representation of $Cl\left( V,Q\right) $, up to
isomorphism, called the \textit{spin representation }of $Cl\left( V,Q\right) 
$. Under the isomorphism%
\begin{equation*}
Cl\left( E,Q_{\left\vert E\right. }\right) \backsimeq Cl\left( V,Q\right) f
\end{equation*}%
the spin representation is the map%
\begin{equation*}
\rho :Cl\left( V,Q\right) \rightarrow End\left( Cl\left( E,Q_{\left\vert
E\right. }\right) \right) 
\end{equation*}%
such that, for all $\varphi \in Cl\left( V,Q\right) $ and $\gamma \in
Cl\left( E,Q_{\left\vert E\right. }\right) $,%
\begin{equation*}
\left( \left( \rho \left( \varphi \right) \right) \left( \gamma \right)
\right) \cdot f=\varphi \cdot \gamma \cdot f\text{;}
\end{equation*}%
its representation space $Cl\left( E,Q_{\left\vert E\right. }\right) $ is
called the \textit{space of spinors} of $\left( V,Q\right) $, denoted by $%
S\left( V,Q\right) $.\smallskip 

\noindent The \textit{space of even }(respectively: \textit{odd}) \textit{%
spinors} of $\left( V,Q\right) $ is%
\begin{equation*}
S\left( V,Q\right) _{+}=Cl\left( E,Q_{\left\vert E\right. }\right) _{+}\text{
(respectively: }S\left( V,Q\right) _{-}=Cl\left( E,Q_{\left\vert E\right.
}\right) _{-}\text{).}
\end{equation*}%
Inside the space of spinors, the subset of \textit{pure spinors} has a very
important geometrical meaning, as we describe in the following.\smallskip

\noindent Let $W$ be a maximal totally isotropic subspace of $V$ and let $%
f_{W}$ be the product of the vectors in a basis of $W$ ($f_{W}$ is well
defined up to a non zero scalar).\smallskip

\noindent It's not hard to show that $Cl\left( V,Q\right) f\cap
f_{W}Cl\left( V,Q\right) $ is a complex vector space of dimension $1$. So we
can pose%
\begin{equation*}
Cl\left( V,Q\right) f\cap f_{W}Cl\left( V,Q\right) =S\left( V,Q\right) _{W}f
\end{equation*}%
where $S(V,Q)_{W}$ denotes a vector subspace of $S(V,Q)$ of dimension $1$%
.\medskip

\noindent \textbf{Definition 2.1 }Any element of $S(V,Q)_{W}\backslash
\left\{ 0\right\} $ is called \textit{representative spinor} of $W$.
Moreover, we call \textit{pure spinor }any element of $S(V,Q)_{W}\backslash
\left\{ 0\right\} $, for some maximal totally isotropic vector subspace $W$
of $V$.\medskip

\noindent It's easy to check that the subset of pure spinors is a projective
variety,\textit{\ }called\textit{\ spinor variety}, and that it is in $1-1$
correspondence with the variety of maximal totally isotropic vector
subspaces of $V$. Furthermore, the spinor variety has two isomorphic
irreducible components, called \textit{even} and\textit{\ odd pure spinors
variety}. From now on we focus our attention on the first one, which we
denote by $S_{h}$.\smallskip

\noindent Let $\mathcal{B}=\left\{ e_{1},...,e_{h},f_{1},...,f_{h}\right\} $
be a basis of $V=E\oplus F$, where $\left\{ e_{1},...,e_{h}\right\} $ is a
basis of $E$ and $\left\{ f_{1},...,f_{h}\right\} $ is a basis of $F$, such
that $B\left( e_{i},f_{j}\right) =\dfrac{\delta _{ij}}{2}$, for all $i,j\in
\left\{ 1,...,h\right\} $. We remark that the matrix $\mathfrak{B}$ of the
form $B$ with respect to $\mathcal{B}$ is%
\begin{equation*}
\mathfrak{B}=%
\begin{bmatrix}
O_{h} & \dfrac{1}{2}I_{h} \\ 
\dfrac{1}{2}I_{h} & O_{h}%
\end{bmatrix}%
\end{equation*}%
where $O_{h}$ and $I_{h}$ are the null matrix and the identity matrix of
size $h$, respectively. Moreover, we pose $f=f_{1}\cdot ...\cdot f_{h}$.

\noindent Let $W$ be a vector subspace of $V$ such that $\dim _{\mathbb{C}%
}W=h$, i.e. $W\in Gr\left( h,2h\right) $, the usual Grassmannian. Thus, we
can associate to $W$ the $h$ by $2h$ matrix%
\begin{equation*}
P=\left[ C_{W}\left\vert D_{W}\right. \right]
\end{equation*}%
where $C_{W},D_{W}\in M\left( h,\mathbb{C}\right) $. In particular, if $%
C_{W} $ is invertible, then we can assume that%
\begin{equation*}
P=\left[ I_{h}\left\vert U_{W}\right. \right]
\end{equation*}%
where $U_{W}=C_{W}^{-1}D_{W}$. So, we have that $W$ is totally isotropic if
and only if%
\begin{equation*}
P\cdot \mathfrak{B}\cdot P^{t}=O_{h}\text{,}
\end{equation*}%
in other words if and only if%
\begin{equation*}
U_{W}=-U_{W}^{t}\text{.}
\end{equation*}%
\noindent We immediately get the following:\medskip

\noindent \textbf{Theorem 2.1} \textit{The generic element of} $S_{h}$ 
\textit{can be represented in block matrix} \textit{form as} $\left[
I_{h}\left\vert U\right. \right] $, \textit{wher}e $U\in M\left( h,\mathbb{C}%
\right) $ \textit{is skew symmetric}.\medskip

\noindent Now, let $U=\left\{ u_{ij}\right\} $ be a skew symmetric matrix of
size $h$ with complex entries and let%
\begin{equation*}
s\left( U\right) =\left( e_{1}+\dsum\limits_{j=1}^{h}u_{1j}f_{j}\right)
\cdot \left( e_{2}+\dsum\limits_{j=1}^{h}u_{2j}f_{j}\right) \cdot ...\cdot
\left( e_{h}+\dsum\limits_{j=1}^{h}u_{hj}f_{j}\right)
\end{equation*}%
be an element of $S_{h}$ in a neighborhood of%
\begin{equation*}
s_{0}=e_{1}\cdot ...\cdot e_{h}\text{.}
\end{equation*}%
We remark that $s\left( U\right) $ and $s_{0}$ are representative spinors of%
\begin{equation*}
W\left( U\right) =\left\langle
e_{1}+\dsum\limits_{j=1}^{h}u_{1j}f_{j},e_{2}+\dsum%
\limits_{j=1}^{h}u_{2j}f_{j},...,e_{h}+\dsum\limits_{j=1}^{h}u_{hj}f_{j}%
\right\rangle
\end{equation*}%
and of $E=W\left( O_{h}\right) $ respectively. By computing $s\left(
U\right) f$ we get the following formula, $\cite{Ang}$ and $\cite{M}$:%
\begin{equation*}
s\left( U\right) =\dsum\limits_{K}Pf_{K}(U)e_{K^{c}}
\end{equation*}%
where $K$ denotes any sequence of integers between $1$ and $h$ of even
lenght, $K^{c}=\left\{ 1,...,h\right\} \backslash K$, $Pf_{K}(U)$ is the
Pfaffian of the submatrix of $U$ made up by rows and columns indexed by $K$,
and $e_{K^{c}}$ is the Clifford product of the $e_{i}$'s, $i\in K^{c}$%
.\smallskip

\noindent In this way we get one of the main tools for our
investigations:\medskip

\noindent \textbf{Theorem 2.2} \textit{All the principal sub-Pfaffians of a
generic skew symmetric matrix of size} $h$ \textit{parametrize a generic
element of }$S_{h}$ \textit{in} $\mathbb{P}^{2^{h-1}-1}\left( \mathbb{C}%
\right) $.\medskip

\noindent Before closing this section we remark that, given%
\begin{equation*}
g=\left[ 
\begin{array}{cc}
g_{11} & g_{12} \\ 
g_{21} & g_{22}%
\end{array}%
\right] \in SO\left( 2h,Q\right)
\end{equation*}%
where $g_{ij}\in M\left( h,\mathbb{C}\right) $, $i,j\in \left\{ 1,2\right\} $
and%
\begin{equation*}
P=\left[ I_{h}\left\vert U\right. \right] \in S_{h}
\end{equation*}%
where $U\in M\left( h,\mathbb{C}\right) $ is skew symmetric, $g$ acts on $P$
as follows:%
\begin{equation*}
g\left( P\right) =\left[ I_{h}\left\vert \left(
g_{11}^{t}+U_{h}g_{12}^{t}\right) ^{-1}\left(
g_{21}^{t}+U_{h}g_{22}^{t}\right) \right. \right] \text{,}
\end{equation*}%
when $\left( g_{11}^{t}+U_{h}g_{12}^{t}\right) ^{-1}$ is defined. As we can
see in $\cite{M}$, this action is \textit{generically} $3$\textit{-transitive%
}, i.e. $Spin\left( 2h,Q\right) $ has an open orbit in $S_{h}\times
S_{h}\times S_{h}$. In order to prove theorem 1.1 part (ii), in section $5$
we provide a proof of this statement based on a new argument: namely we
consider $3$ points of $S_{h}$ that are in the same parametrization (see
theorem 5.1).

\section{Higher secant varieties$\protect\medskip $}

\noindent Let $X$ $\subseteq \mathbb{P}^{N}\left( \mathbb{C}\right) $ be a $%
d $-dimensional projective variety.\smallskip

\noindent We pose the following:\medskip

\noindent \textbf{Definition 3.1 }The $k$\textit{-secant variety }$\sigma
_{k}\left( X\right) $ is the Zariski closure of the union of all linear
spaces spanned by $k$ points of $X$, that is%
\begin{equation*}
\sigma _{k}\left( X\right) =\overline{\underset{x_{1},...,x_{k}\in X}{\cup }%
\left\langle x_{1},...,x_{k}\right\rangle }\text{.}
\end{equation*}%
If $X$ $\subseteq \mathbb{P}^{N}\left( \mathbb{C}\right) $ is
non-degenerate, i.e. is not contained in any hyperplane, then we have the
following estimate on the dimension of $\sigma _{k}\left( X\right) $:%
\begin{equation*}
\dim _{\mathbb{C}}\sigma _{k}\left( X\right) \leq \min \left\{
kd+k-1,N\right\} \text{.}
\end{equation*}%
The problem of determining when the dimension of the secant variety $\sigma
_{k}\left( X\right) $ reaches this upper bound is called \textit{defectivity
problem }for $X$. In this sense we have the following:\medskip

\noindent \textbf{Definition 3.2 }Let $X$ $\subseteq \mathbb{P}^{N}\left( 
\mathbb{C}\right) $ be a non-degenerate projective variety of dimension $d$%
.\smallskip

\noindent 1. If $\dim _{\mathbb{C}}\sigma _{k}\left( X\right) =\min \left\{
kd+k-1,N\right\} $ then we say that $\sigma _{k}\left( X\right) $ has the 
\textit{expected dimension}.\smallskip

\noindent 2. If $\dim _{\mathbb{C}}\sigma _{k}\left( X\right) <\min \left\{
kd+k-1,N\right\} $ then we say that $X$ has a \textit{defective }$k$\textit{%
-secant variety} and that%
\begin{equation*}
\delta _{k}=\min \left\{ k\dim _{\mathbb{C}}X+k-1,N\right\} -\dim _{\mathbb{C%
}}\sigma _{k}\left( X\right)
\end{equation*}%
is its $k$-\textit{defect}.\smallskip

\noindent 3. If there's a $k$ such that $X$ is $k$-defective then we say
that $X$ is \textit{defective}.\textit{\medskip }

\noindent Now we recall the main tool to compute the dimensions of higher
secant varieties:\medskip

\noindent \textbf{Lemma 3.1} (\textbf{Terracini}, 1911) \textit{Let} $%
X\subset \mathbb{P}^{N}\left( \mathbb{C}\right) $ \textit{be a projective
variety and let} $z$ \textit{be a generic smooth point of }$\sigma
_{k}\left( X\right) $. \textit{Then the projective tangent space to }$\sigma
_{k}\left( X\right) $\textit{\ at }$z$ \textit{is given by}%
\begin{equation*}
\widetilde{T}_{z}\sigma _{k}\left( X\right) =\left\langle \widetilde{T}%
_{x_{1}}X,...,\widetilde{T}_{x_{k}}X\right\rangle
\end{equation*}%
\textit{where} $x_{1},...,x_{k}$ \textit{are generic smooth points of }$X$ 
\textit{such that }$z\mathbb{\in }\left\langle x_{1},...,x_{k}\right\rangle $
\textit{and} $\widetilde{T}_{x_{i}}X$ \textit{denotes} \textit{the
projective tangent space to }$X$\textit{\ at }$x_{i}$.\medskip

\noindent By upper semicontinuity, we immediately get an argument to prove
that a variety isn't defective:\medskip

\noindent \textbf{Corollary 3.2 }\textit{Let }$x_{1},...,x_{k}\in X$\textit{%
\ be smooth points such that }$\widetilde{T}_{x_{1}}X,...,\widetilde{T}%
_{x_{k}}X$ \textit{are linearly independent}, \textit{or else}%
\begin{equation*}
\left\langle \widetilde{T}_{x_{1}}X,...,\widetilde{T}_{x_{k}}X\right\rangle =%
\mathbb{P}^{N}\left( \mathbb{C}\right) \text{.}
\end{equation*}%
\textit{Then} $\sigma _{k}\left( X\right) $ \textit{has the expected
dimension}.\medskip

\noindent \textit{Terracini's lemma} also provides a method to show that $X$
has a defective $k$-secant variety. More precisely, we have the
following:\medskip

\noindent \textbf{Corollary 3.3 (}$\cite{C-C}$\textbf{) }\textit{Let }$%
d=\dim _{\mathbb{C}}X$ \textit{and let us suppose that}%
\begin{equation}
kd+k-1\leq N\text{.}  \label{7}
\end{equation}%
\noindent \textit{If there exists a rational normal curve of }$X$, \textit{%
embedded in }$\mathbb{P}^{2k-2}\left( \mathbb{C}\right) $ \textit{and
containing }$k$ \textit{general points of }$X$, \textit{then }$\sigma
_{k}\left( X\right) $\textit{\ hasn't the expected dimension}.\medskip

\noindent \textbf{Proof. }Let $x_{1},...,x_{k}$ be general points of $X$
satisfying the hypothesis and let $T_{x_{1}}X,...,T_{x_{k}}X$ be the affine
tangent spaces at such points. We get that, for all $i\in \left\{
1,...,k\right\} $,%
\begin{equation}
\dim _{%
%TCIMACRO{\U{2102} }%
%BeginExpansion
\mathbb{C}
%EndExpansion
}\left( T_{x_{i}}X\cap 
%TCIMACRO{\U{2102} }%
%BeginExpansion
\mathbb{C}
%EndExpansion
^{2k-1}\right) =2  \label{8}
\end{equation}%
\noindent because $T_{x_{i}}X$ contains the affine tangent space to the
curve at $x_{i}$. Now, let $\pi _{|\left\langle
T_{x_{1}}X,...,T_{x_{k}}X\right\rangle }$ be the restriction to $%
\left\langle T_{x_{1}}X,...,T_{x_{k}}X\right\rangle $ of the canonical
projection%
\begin{equation*}
\pi :%
%TCIMACRO{\U{2102} }%
%BeginExpansion
\mathbb{C}
%EndExpansion
^{N+1}\rightarrow 
%TCIMACRO{\U{2102} }%
%BeginExpansion
\mathbb{C}
%EndExpansion
^{N+1}/%
%TCIMACRO{\U{2102} }%
%BeginExpansion
\mathbb{C}
%EndExpansion
^{2k-1}\text{.}
\end{equation*}%
We remark that $\pi $ is a linear mapping between vector spaces, thus%
\begin{eqnarray*}
\dim _{%
%TCIMACRO{\U{2102} }%
%BeginExpansion
\mathbb{C}
%EndExpansion
}\left\langle T_{x_{1}}X,...,T_{x_{k}}X\right\rangle &=&\dim _{%
%TCIMACRO{\U{2102} }%
%BeginExpansion
\mathbb{C}
%EndExpansion
}\ker \pi _{|\left\langle T_{x_{1}}X,...,T_{x_{k}}X\right\rangle }+\smallskip
\\
&&+\dim _{%
%TCIMACRO{\U{2102} }%
%BeginExpansion
\mathbb{C}
%EndExpansion
}\func{Im}\pi _{|\left\langle T_{x_{1}}X,...,T_{x_{k}}X\right\rangle
\smallskip } \\
&=&\dim _{%
%TCIMACRO{\U{2102} }%
%BeginExpansion
\mathbb{C}
%EndExpansion
}\left( \left\langle T_{x_{1}}X,...,T_{x_{k}}X\right\rangle \cap 
%TCIMACRO{\U{2102} }%
%BeginExpansion
\mathbb{C}
%EndExpansion
^{2k-1}\right) +\smallskip \\
&&+\dim _{%
%TCIMACRO{\U{2102} }%
%BeginExpansion
\mathbb{C}
%EndExpansion
}\left\langle \pi \left( T_{x_{1}}X\right) ,...,\pi \left( T_{x_{k}}X\right)
\right\rangle \text{.}
\end{eqnarray*}%
and then%
\begin{equation}
\dim _{%
%TCIMACRO{\U{2102} }%
%BeginExpansion
\mathbb{C}
%EndExpansion
}\left\langle T_{x_{1}}X,...,T_{x_{k}}X\right\rangle \leq
2k-1+\dsum\limits_{i=1}^{k}\dim _{%
%TCIMACRO{\U{2102} }%
%BeginExpansion
\mathbb{C}
%EndExpansion
}\pi \left( T_{x_{i}}X\right) \text{.}  \label{10}
\end{equation}%
\noindent Now, let $\pi _{|T_{x_{i}}X}$ be the restriction of $\pi $ to $%
T_{x_{i}}X$; from $\left( \ref{10}\right) $ and $\left( \ref{8}\right) $ we
get that%
\begin{eqnarray}
\dim _{%
%TCIMACRO{\U{2102} }%
%BeginExpansion
\mathbb{C}
%EndExpansion
}\left\langle T_{x_{1}}X,...,T_{x_{k}}X\right\rangle &\leq
&2k-1+\dsum\limits_{i=1}^{k}\dim _{%
%TCIMACRO{\U{2102} }%
%BeginExpansion
\mathbb{C}
%EndExpansion
}T_{x_{i}}X-\dim _{%
%TCIMACRO{\U{2102} }%
%BeginExpansion
\mathbb{C}
%EndExpansion
}\left( T_{x_{i}}X\cap 
%TCIMACRO{\U{2102} }%
%BeginExpansion
\mathbb{C}
%EndExpansion
^{2k-1}\right)  \notag \\
&=&2k-1+k\left[ (\dim _{%
%TCIMACRO{\U{2102} }%
%BeginExpansion
\mathbb{C}
%EndExpansion
}X+1)-2\right] \smallskip  \label{11} \\
&=&k\dim _{%
%TCIMACRO{\U{2102} }%
%BeginExpansion
\mathbb{C}
%EndExpansion
}X+k-1\text{.}  \notag
\end{eqnarray}

\noindent Finally, let $\widehat{\sigma _{k}\left( X\right) }$ be the affine
cone over $\sigma _{k}\left( X\right) $; by using $\left( \ref{7}\right) $
we immediately have that the expected dimension for $\widehat{\sigma
_{k}\left( X\right) }$ is%
\begin{equation*}
\exp \dim _{\mathbb{C}}\widehat{\sigma _{k}\left( X\right) }=k\dim _{%
%TCIMACRO{\U{2102} }%
%BeginExpansion
\mathbb{C}
%EndExpansion
}X+k\text{.}
\end{equation*}%
Then, from Terracini's lemma and from $\left( \ref{11}\right) $, we get that%
\begin{equation*}
\dim _{\mathbb{C}}\widehat{\sigma _{k}\left( X\right) }<\exp \dim _{\mathbb{C%
}}\widehat{\sigma _{k}\left( X\right) }\text{,}
\end{equation*}%
i.e. $X$ has a defective $k$-secant variety.\textit{\ }\hfill $\square $

\section{A probabilistic algorithm for the secant defect of spinor varieties 
$\protect\medskip $}

\noindent To deal with our problem, we constructed a probabilistic algorithm
through the \textit{Macaulay2 }computation system, $\cite{M2}$.\smallskip

\noindent The script of the algorithm is given below:\bigskip

\noindent {\footnotesize h = value read "h?"}

{\footnotesize \noindent k = value read "k?"}

{\footnotesize \noindent p = floor(h*(h-1)/2)}

{\footnotesize \noindent R = QQ[x\_0..x\_(p-1)]}

{\footnotesize \noindent X = vars R}

{\footnotesize \noindent M = X -\TEXTsymbol{>} genericSkewMatrix(R,x\_0,h)}

{\footnotesize \noindent par = X -\TEXTsymbol{>} apply(floor(h/2)+1,i-%
\TEXTsymbol{>}generators pfaffians(2*i,M(X)))}

{\footnotesize \noindent f = l -\TEXTsymbol{>} (a=l\#0;for i from 1 to
\#(l)-1 do(a=a\TEXTsymbol{\vert}(l\#i););a)}

{\footnotesize \noindent S = f(par(X))}

{\footnotesize \noindent J = jacobian S}

{\footnotesize \noindent g = l -\TEXTsymbol{>} (a=l\#0;for i from 1 to
\#(l)-1 do(a=a\TEXTsymbol{\vert}\TEXTsymbol{\vert}(l\#i););a)}

{\footnotesize \noindent punti = apply(k,i-\TEXTsymbol{>}for j from 1 to p
list random(1000))}

{\footnotesize \noindent puntibis = apply(k,i-\TEXTsymbol{>}%
matrix\{punti\#i\})}

{\footnotesize \noindent Spunti = apply(k,i-\TEXTsymbol{>}%
substitute(S,matrix(R,\{flatten entries puntibis\#i\})))}

{\footnotesize \noindent Jpunti = apply(k,i-\TEXTsymbol{>}%
substitute(J,matrix(R,\{flatten entries puntibis\#i\})))}

{\footnotesize \noindent JS = apply(k,i-\TEXTsymbol{>}(Spunti\#i)\TEXTsymbol{%
\vert}\TEXTsymbol{\vert}(Jpunti\#i))}

{\footnotesize \noindent JJS = g(JS)}

{\footnotesize \noindent rank JJS.\bigskip }

\noindent This algorithm is based on \textit{Terracini's lemma }and on the
fact that \textit{Pfaffians }parametrize $S_{h}$; moreover it was conceived
for every $h$ and $k$ integers, where $h=\dfrac{1}{2}\dim _{%
%TCIMACRO{\U{2102} }%
%BeginExpansion
\mathbb{C}
%EndExpansion
}V$.\smallskip

\noindent The main steps of our algorithm are the following:\medskip

\noindent 1. \textit{Preliminaries.\smallskip }

\noindent Given $h$, $k$ and further computed the dimension of $S_{h}$%
\begin{equation*}
p=\dfrac{h(h-1)}{2}\text{,}
\end{equation*}%
we define the polinomial ring $R$ with rational coefficients in the
variables $\left\{ x_{0},...,x_{p-1}\right\} $.\medskip

\noindent 2. \textit{Parametrization of }$S_{h}$\textit{.\smallskip }

\noindent In order to parametrize the variety of even pure spinors, we
construct the function%
\begin{equation*}
M:\mathcal{M}_{\left( 1,p\right) }\left( 
%TCIMACRO{\U{211a} }%
%BeginExpansion
\mathbb{Q}
%EndExpansion
\right) \rightarrow \mathcal{M}_{\left( h,h\right) }\left( 
%TCIMACRO{\U{211a} }%
%BeginExpansion
\mathbb{Q}
%EndExpansion
\right)
\end{equation*}%
defined by%
\begin{equation*}
X=\left( x_{0},...,x_{p-1}\right) \rightarrow M\left( X\right) =\left( 
\begin{array}{ccccc}
0 & x_{0} & x_{1} & \cdots & x_{h-1} \\ 
-x_{0} & 0 & \ddots & \ddots & \vdots \\ 
-x_{1} & \ddots & \ddots & \ddots & \vdots \\ 
\vdots & \ddots & \ddots & \ddots & x_{p-1} \\ 
-x_{h-1} & \cdots & \cdots & -x_{p-1} & 0%
\end{array}%
\right)
\end{equation*}%
and then we compute the principal sub-Pfaffians\textit{\ }of this matrix by
using the function%
\begin{equation*}
par:\mathcal{M}_{\left( 1,p\right) }\left( 
%TCIMACRO{\U{211a} }%
%BeginExpansion
\mathbb{Q}
%EndExpansion
\right) \rightarrow \mathcal{M}_{\left( 1,2^{h-1}\right) }\left( 
%TCIMACRO{\U{211a} }%
%BeginExpansion
\mathbb{Q}
%EndExpansion
\right)
\end{equation*}%
such that%
\begin{equation*}
X=\left( x_{0},...,x_{p-1}\right) \rightarrow par\left( X\right) =\left( 
\text{principal sub-Pfaffians\textit{\ }of }M\left( X\right) \right) \text{.}
\end{equation*}

\noindent 3. \textit{Definition of }$S_{h}$\textit{.\smallskip }

\noindent From the theorem 2.2 we obtain that $S_{h}$ is the image of the
function \textit{par}, i.e. it belongs to $\mathcal{M}_{\left(
1,2^{h-1}\right) }\left( 
%TCIMACRO{\U{211a} }%
%BeginExpansion
\mathbb{Q}
%EndExpansion
\right) $:%
\begin{equation*}
S=par\left( X\right) =\left( s_{i}\right) _{i=0,...,2^{h-1}-1}\text{.}
\end{equation*}%
We observe that \textit{par}, being defined through \textit{apply}, produces
a list of $\left\lfloor \dfrac{h}{2}\right\rfloor +1$ row matrices; by means
of the function $f$ we juxtapose all Pfaffians in one row matrix.\medskip

\noindent 4. \textit{Computation of the jacobian matrix of the
parametrization.\smallskip }

\noindent Applying \textit{jacobian} to $S$ we get the following $p$ by $%
2^{h-1}$ matrix:%
\begin{equation*}
J=\left( \partial _{j}s_{i}\right) _{i=0,...,2^{h-1}-1;j=0,...,p-1}\text{.}
\end{equation*}

\noindent 5. \textit{Choice of }$k$ \textit{random points in }$S_{h}$ 
\textit{and computation of their coordinates.\smallskip }

\noindent In order to study $\sigma _{k}\left( S_{h}\right) $, we have to
choose $k$ elements of $S_{h}$: so, we consider a list of $k$ sets (\textit{%
punti}) of $p$ random rational numbers and we construct the corresponding
skew symmetric $h$ by $h$ matrices; then we compute the principal
sub-Pfaffians of these matrices\textit{. }In this way we get a list (\textit{%
Spunti}) composed of the parametric coordinates of the $k$ selected points%
\textit{.}%
\begin{equation*}
\begin{array}{lll}
punti & = & \left\{ punti_{0},...,punti_{k-1}\right\} \smallskip \\ 
punti_{i} & = & \left( q_{0}^{i},...,q_{p-1}^{i}\right) \text{, }%
q_{j}^{i}\in 
%TCIMACRO{\U{211a} }%
%BeginExpansion
\mathbb{Q}
%EndExpansion
\text{ random, }q_{j}^{i}\leq 1000\smallskip \\ 
Spunti & = & \left\{ S\left( punti_{0}\right) ,...,S\left(
punti_{k-1}\right) \right\} =\left\{ P_{0},...,P_{k-1}\right\}%
\end{array}%
\end{equation*}

\noindent 6. \textit{Construction of the affine tangent spaces to }$S_{h}$ 
\textit{at the }$k$ \textit{points.\smallskip }

\noindent Now we evaluate the jacobian matrix $J$ at the points under
consideration. Thus we obtain a list (\textit{Jpunti}) of matrices whose
images correspond to the vector tangent spaces to $S_{h}$; placing the row
made up of the coordinates of one point before the corresponding jacobian
matrix we get the affine tangent space to $S_{h}$ at such point.%
\begin{equation*}
\begin{array}{lll}
Jpunti & = & \left\{ J|_{X=punti_{0}},...,J|_{X=punti_{k-1}}\right\}
=\left\{ J_{0},...,J_{k-1}\right\} \smallskip \\ 
JS & = & \left\{ P_{0}|J_{0},...,P_{k-1}|J_{k-1}\right\} =\left\{
JS_{0},...,JS_{k-1}\right\}%
\end{array}%
\end{equation*}

\noindent 7. \textit{Computation of the dimension of }$\sigma _{k}\left(
S_{h}\right) $\textit{.\smallskip }

\noindent Finally, we arrange in columns the $\left( p+1\right) $ by $%
2^{h-1} $ matrices $JS_{0},...,JS_{k-1}$ and we obtain the $k\left(
p+1\right) $ by $2^{h-1}$ matrix $JJS$ associated with the span of the
affine tangent spaces. From \textit{Terracini's Lemma} we get that the rank
of $JJS$ produces the affine dimension of $\sigma _{k}\left( S_{h}\right) $;
subtracting $1$ to the output we get the required dimension.%
\begin{eqnarray*}
&&%
\begin{array}{lllll}
g & : & \left\{ \text{lists of matrices}\right\} & \rightarrow & \left\{ 
\text{matrices}\right\} \smallskip \\ 
&  & B=\left\{ B_{1},B_{2},...\right\} & \rightarrow & g\left( B\right)
=\left( B_{1}|B_{2}|...\right) ^{t}\smallskip \\ 
g\left( JS\right) & = & \left( 
\begin{array}{c}
JS_{0} \\ 
\vdots \\ 
\vdots \\ 
JS_{k-1}%
\end{array}%
\right) & = & JJS%
\end{array}
\\
&&\text{OUTPUT }rank\left( JJS\right)
\end{eqnarray*}

\noindent \textbf{Remark 4.1 }If the achieved value coincides with the
expected dimension of $\sigma _{k}\left( S_{h}\right) $, i.e. if $JJS$ has
maximum rank, then we can be sure that the actual dimension is that value
(corollary 3.2); otherwise\ we need other checks to say that $S_{h}$ is $k$%
-defective.

\noindent Thus we can say that our algorithm is probabilistic.\medskip

\noindent It's not hard to check, by direct computations, that, if $h\leq 5$%
, then $S_{h}$ isn't defective, $\cite{Ang}$ and $\cite{K}$. So we used this
algorithm from the stage $(h,k)=\left( 6,2\right) $ to the stage $%
(h,k)=\left( 9,5\right) $: beyond these values the memory of the computer
was used up.\medskip

\noindent Our results are summarized as follows.

\begin{equation*}
\fbox{$\mathbf{k=2}$}
\end{equation*}%
\begin{equation*}
\begin{array}{cccccc}
h & p & N & \exp \dim \sigma _{k}\left( S_{h}\right) & \dim \sigma
_{k}\left( S_{h}\right) & \text{defective} \\ 
6 & 15 & 31 & 31 & 31 & \text{NO} \\ 
7 & 21 & 63 & 43 & 43 & \text{NO} \\ 
8 & 28 & 127 & 57 & 57 & \text{NO} \\ 
9 & 36 & 255 & 73 & 73 & \text{NO} \\ 
10 & 45 & 511 & 91 & 91 & \text{NO} \\ 
11 & 55 & 1023 & 111 & 111 & \text{NO}%
\end{array}%
\end{equation*}%
\begin{equation*}
\fbox{$\mathbf{k=3}$}
\end{equation*}%
\begin{equation*}
\begin{array}{cccccc}
h & p & N & \exp \dim \sigma _{k}\left( S_{h}\right) & \dim \sigma
_{k}\left( S_{h}\right) & \text{defective} \\ 
7 & 21 & 63 & 63 & 58 & \text{YES\QQfnmark{%
see theorem 5.5.}} \\ 
8 & 28 & 127 & 86 & 85 & \text{YES\QQfnmark{%
see theorem 5.3.}} \\ 
9 & 36 & 255 & 110 & 110 & \text{NO} \\ 
10 & 45 & 511 & 137 & 137 & \text{NO} \\ 
11 & 55 & 1023 & 167 & 167 & \text{NO} \\ 
12 & 66 & 2047 & 200 & 200 & \text{NO}%
\end{array}%
\QQfntext{-1}{
see theorem 5.5.}
\QQfntext{1}{
see theorem 5.3.}
\end{equation*}%
\begin{equation*}
\fbox{$\mathbf{k=4}$}
\end{equation*}%
\begin{equation*}
\begin{array}{cccccc}
h & p & N & \exp \dim \sigma _{k}\left( S_{h}\right) & \dim \sigma
_{k}\left( S_{h}\right) & \text{defective} \\ 
7 & 21 & 63 & 63 & 63 & \text{NO} \\ 
8 & 28 & 127 & 115 & 111 & \text{YES\QQfnmark{%
see corollary 5.4.}} \\ 
9 & 36 & 255 & 147 & 147 & \text{NO} \\ 
10 & 45 & 511 & 183 & 183 & \text{NO}%
\end{array}%
\QQfntext{0}{
see corollary 5.4.}
\end{equation*}%
\begin{equation*}
\fbox{$\mathbf{k=5}$}
\end{equation*}%
\begin{equation*}
\begin{array}{cccccc}
h & p & N & \exp \dim \sigma _{k}\left( S_{h}\right) & \dim \sigma
_{k}\left( S_{h}\right) & \text{defective} \\ 
8 & 28 & 127 & 127 & 127 & \text{NO} \\ 
9 & 36 & 255 & 184 & 184 & \text{NO}%
\end{array}%
\end{equation*}%
\medskip

\noindent The last three tables provide a proof of theorem 1.1 part (i) till 
$h=12$ and even some cases more.\smallskip

\noindent In the first table we can see that, if $6\leq h\leq 11$, then $%
\sigma _{2}\left( S_{h}\right) $ has the expected dimension; this fact
agrees with already known theoretical results, $\cite{K}$.\smallskip

\noindent However, we found some "anomalies" when $(h,k)\in \left\{ \left(
7,3\right) ,\left( 8,3\right) ,\left( 8,4\right) \right\} $. So, we supposed
that actually these varieties haven't the expected dimension. Indeed, in the
next section we explain, from a theoretical point of view, that $S_{8}$ has
a defective $3$-secant variety and a defective $4$-secant variety and that $%
S_{7}$ has a defective $3$-secant variety. Hence we get a proof of theorem
1.1 part (ii).

\section{The defective cases\emph{\protect\medskip }}

\noindent In order to prove that $\sigma _{3}\left( S_{8}\right) $ and $%
\sigma _{4}\left( S_{8}\right) $ haven't the expected dimension, we proceed
as follows.\medskip

\noindent Let assume that $h$ is an even number, $h=2m$. With the notations
of section $2$, let%
\begin{equation*}
s_{0}=e_{1}\cdot ...\cdot e_{h}\text{, }s_{1}=\dprod\limits_{i=1}^{m}\left(
1+e_{2i-1}\cdot e_{2i}\right) \text{, }s_{2}=\dprod\limits_{i=1}^{m}\left(
1-e_{2i-1}\cdot e_{2i}\right)
\end{equation*}%
be elements of $S_{h}$: they are \textit{representative spinors} of the
maximal totally isotropic subspaces%
\begin{eqnarray*}
E &=&\left\langle e_{1},...,e_{h}\right\rangle \smallskip \\
G &=&\left\langle e_{2i-1}+f_{2i},e_{2i}-f_{2i-1},1\leq i\leq m\right\rangle
\smallskip \\
H &=&\left\langle e_{2i-1}-f_{2i},e_{2i}+f_{2i-1},1\leq i\leq m\right\rangle
\end{eqnarray*}%
respectively. Their corresponding $h$ by $2h$ matrices are%
\begin{eqnarray*}
P_{0} &=&\left[ I_{h}\left\vert O_{h}\right. \right] \smallskip \\
P_{1} &=&\left[ I_{h}\left\vert J_{m}\right. \right] \smallskip \\
P_{2} &=&\left[ I_{h}\left\vert -J_{m}\right. \right]
\end{eqnarray*}%
where $J_{m}$ denotes the skew symmetric matrix of size $h$ made up of $m$
diagonal blocks like $%
\begin{pmatrix}
0 & 1 \\ 
-1 & 0%
\end{pmatrix}%
$.\medskip

\noindent \textbf{Theorem 5.1 }\textit{The orbit of} $\left(
P_{0},P_{1},P_{2}\right) $ \textit{is open in} $S_{h}\times S_{h}\times
S_{h} $.\medskip

\noindent \textbf{Proof. }Let consider the function%
\begin{equation*}
\begin{array}{ccccc}
f & : & SO(2h,Q) & \rightarrow & S_{h}\times S_{h}\times S_{h}\smallskip \\ 
&  & g & \rightarrow & \left( g\left( P_{0}\right) ,g\left( P_{1}\right)
,g\left( P_{2}\right) \right)%
\end{array}%
\end{equation*}%
where%
\begin{eqnarray*}
g\left( P_{0}\right) &=&\left[ I_{h}\left\vert \left( g_{11}^{t}\right)
^{-1}g_{21}^{t}\right. \right] \smallskip \\
g\left( P_{1}\right) &=&\left[ I_{h}\left\vert \left(
g_{11}^{t}+J_{m}g_{12}^{t}\right) ^{-1}\left(
g_{21}^{t}+J_{m}g_{22}^{t}\right) \right. \right] \smallskip \\
g\left( P_{2}\right) &=&\left[ I_{h}\left\vert \left(
g_{11}^{t}-J_{m}g_{12}^{t}\right) ^{-1}\left(
g_{21}^{t}-J_{m}g_{22}^{t}\right) \right. \right] \text{;}
\end{eqnarray*}%
we remark that%
\begin{equation*}
\func{Im}f=\left\{ \left( g\left( P_{0}\right) ,g\left( P_{1}\right)
,g\left( P_{2}\right) \right) \left\vert \text{ }g\in SO(2h,Q)\right.
\right\}
\end{equation*}%
is the orbit of $\left( P_{0},P_{1},P_{2}\right) $. Taking $g=I_{2h}$, the
tangent map of $f$ at the point $g$ is:%
\begin{equation*}
df_{I_{2h}}:so\left( 2h,Q\right) \rightarrow T_{\left(
P_{0},P_{1},P_{2}\right) }\left[ S_{h}\times S_{h}\times S_{h}\right] \text{,%
}
\end{equation*}%
where $so\left( 2h,Q\right) $ is the Lie algebra of $SO(2h,Q)$, that is:%
\begin{equation*}
so\left( 2h,Q\right) =\left\{ A\in SO(2h,Q)\left\vert \text{ }A^{t}\mathfrak{%
B}=-\mathfrak{B}A\right. \right\} \text{.}
\end{equation*}%
We have that $\func{Im}df_{I_{2h}}$ is the tangent space to the orbit of $%
\left( P_{0},P_{1},P_{2}\right) $ at $\left( P_{0},P_{1},P_{2}\right) $. Our
aim is to show that $df_{I_{2h}}$ is surjective, or that%
\begin{eqnarray*}
\dim _{%
%TCIMACRO{\U{2102} }%
%BeginExpansion
\mathbb{C}
%EndExpansion
}\ker df_{I_{2h}} &=&\dim _{%
%TCIMACRO{\U{2102} }%
%BeginExpansion
\mathbb{C}
%EndExpansion
}so\left( 2h,Q\right) -\dim _{%
%TCIMACRO{\U{2102} }%
%BeginExpansion
\mathbb{C}
%EndExpansion
}\func{Im}df_{I_{2h}}\smallskip \\
&=&\frac{2h(2h-1)}{2}-\frac{3h(h-1)}{2}\smallskip \\
&=&\frac{h(h+1)}{2}\text{.}
\end{eqnarray*}%
In order to study $\ker df_{I_{2h}}$, we use the first-order Taylor
expansion of $f=\left( f_{1},f_{2},f_{3}\right) $ about $I_{2h}$. So, let $%
H\in so\left( 2h,Q\right) $, i.e.%
\begin{equation*}
H=%
\begin{bmatrix}
H_{11} & H_{12} \\ 
H_{21} & H_{22}%
\end{bmatrix}%
\end{equation*}%
with $H_{ij}\in M\left( h,%
%TCIMACRO{\U{2102} }%
%BeginExpansion
\mathbb{C}
%EndExpansion
\right) $, $i,j\in \left\{ 1,2\right\} $, such that $H_{11}^{t}=-H_{22}$ and 
$H_{12}$, $H_{21}$ are skew symmetric; we get that%
\begin{equation*}
\begin{array}{lll}
f_{1}\left( I_{2h}+H\right) & = & \left[ I_{h}\left\vert \left(
I_{h}+H_{11}^{t}\right) ^{-1}H_{21}^{t}\right. \right] =\left[
I_{h}\left\vert H_{21}^{t}+...\right. \right] \medskip \\ 
f_{2}\left( I_{2h}+H\right) & = & \left[ I_{h}\left\vert \left(
I_{h}+H_{11}^{t}+J_{m}H_{12}^{t}\right) ^{-1}\left( H_{21}^{t}+J_{m}\left(
I_{h}+H_{22}^{t}\right) \right) \right. \right] \smallskip \\ 
& = & \left[ I_{h}\left\vert
J_{m}+H_{21}^{t}+J_{m}H_{22}^{t}-H_{11}^{t}J_{m}-J_{m}H_{12}^{t}J_{m}+...%
\right. \right] \medskip \\ 
f_{3}\left( I_{2h}+H\right) & = & \left[ I_{h}\left\vert \left(
I_{h}+H_{11}^{t}-J_{m}H_{12}^{t}\right) ^{-1}\left( H_{21}^{t}-J_{m}\left(
I_{h}+H_{22}^{t}\right) \right) \right. \right] \smallskip \\ 
& = & \left[ I_{h}\left\vert
-J_{m}+H_{21}^{t}-J_{m}H_{22}^{t}+H_{11}^{t}J_{m}-J_{m}H_{12}^{t}J_{m}+...%
\right. \right]%
\end{array}%
\end{equation*}%
and then we have that%
\begin{equation*}
\ker df_{I_{2h}}=\left\{ H\in so\left( 2h,Q\right) \left\vert \text{ }%
\begin{array}{l}
H_{21}^{t}=0 \\ 
H_{21}^{t}+J_{m}H_{22}^{t}-H_{11}^{t}J_{m}-J_{m}H_{12}^{t}J_{m}=0 \\ 
H_{21}^{t}-J_{m}H_{22}^{t}+H_{11}^{t}J_{m}-J_{m}H_{12}^{t}J_{m}=0%
\end{array}%
\right. \right\} \text{.}
\end{equation*}%
A direct computation $\left( \cite{Ang}\right) $ shows that%
\begin{equation*}
\ker df_{I_{2h}}=\left\{ H\in so\left( 2h,Q\right) \left\vert \text{ }%
\begin{array}{l}
H_{21}^{t}=H_{12}^{t}=0 \\ 
J_{m}H_{22}^{t}=\left( J_{m}H_{22}^{t}\right) ^{t}%
\end{array}%
\right. \right\} \text{,}
\end{equation*}%
thus%
\begin{equation*}
\dim _{%
%TCIMACRO{\U{2102} }%
%BeginExpansion
\mathbb{C}
%EndExpansion
}\ker df_{I_{2h}}=\dim _{%
%TCIMACRO{\U{2102} }%
%BeginExpansion
\mathbb{C}
%EndExpansion
}\left\{ A\in M\left( h,%
%TCIMACRO{\U{2102} }%
%BeginExpansion
\mathbb{C}
%EndExpansion
\right) \left\vert \text{ }J_{m}A=\left( J_{m}A\right) ^{t}\right. \right\} 
\text{.}
\end{equation*}%
Now, by using induction on $m$, where $m=\dfrac{h}{2}$ and $m\geq 1$, we can
prove that%
\begin{equation*}
\dim _{%
%TCIMACRO{\U{2102} }%
%BeginExpansion
\mathbb{C}
%EndExpansion
}\left\{ A\in M\left( h,%
%TCIMACRO{\U{2102} }%
%BeginExpansion
\mathbb{C}
%EndExpansion
\right) \left\vert \text{ }J_{m}A=\left( J_{m}A\right) ^{t}\right. \right\} =%
\frac{h(h+1)}{2}\text{.}
\end{equation*}%
It's not difficult to check the statement for $m=1$.

\noindent Hence, assume the result to be proved till $m$, we want to show
that it holds also for $m+1$. We remark that every $A\in M\left( h+2,%
%TCIMACRO{\U{2102} }%
%BeginExpansion
\mathbb{C}
%EndExpansion
\right) $ can be written as%
\begin{equation*}
A=%
\begin{bmatrix}
\mathcal{A} & 
\begin{array}{c}
\mathcal{B}_{1} \\ 
\vdots \\ 
\mathcal{B}_{m}%
\end{array}
\\ 
\begin{array}{ccc}
\mathcal{C}_{1} & \ldots & \mathcal{C}_{m}%
\end{array}
& \mathcal{D}%
\end{bmatrix}%
\end{equation*}%
with $\mathcal{A\in }M\left( h,%
%TCIMACRO{\U{2102} }%
%BeginExpansion
\mathbb{C}
%EndExpansion
\right) $ and $\mathcal{B}_{1},...,\mathcal{B}_{m},\mathcal{C}_{1},...,%
\mathcal{C}_{m},\mathcal{D\in }M\left( 2,%
%TCIMACRO{\U{2102} }%
%BeginExpansion
\mathbb{C}
%EndExpansion
\right) $. Thus, by the inductive hypothesis we get that%
\begin{equation*}
\dim _{%
%TCIMACRO{\U{2102} }%
%BeginExpansion
\mathbb{C}
%EndExpansion
}\left\{ A\in M\left( h+2,%
%TCIMACRO{\U{2102} }%
%BeginExpansion
\mathbb{C}
%EndExpansion
\right) \left\vert \text{ }J_{m+1}A=\left( J_{m+1}A\right) ^{t}\right.
\right\} =\frac{\left( h+2\right) \left( h+3\right) }{2}
\end{equation*}%
which concludes the proof. \hfill $\square $\medskip

\noindent \textbf{Corollary} \textbf{5.2} \textit{If }$h=2m$\textit{\ then}%
\begin{equation*}
s_{0}=e_{1}\cdot ...\cdot e_{h}\text{, }s_{1}=\dprod\limits_{i=1}^{m}\left(
1+e_{2i-1}\cdot e_{2i}\right) \text{, }s_{2}=\dprod\limits_{i=1}^{m}\left(
1-e_{2i-1}\cdot e_{2i}\right)
\end{equation*}%
\textit{are general points of} $S_{h}$.\bigskip

\noindent Now we are ready to prove the following:\medskip

\noindent \textbf{Theorem} \textbf{5.3} \textit{The} \textit{variety }$S_{8}$
\textit{is} $3$\textit{-defective and }$\delta _{3}=1$.\medskip

\noindent \textbf{Proof. }From corollary 5.2 we get that%
\begin{equation*}
s_{0}=e_{1}\cdot ...\cdot e_{8}\text{, }s_{1}=\dprod\limits_{i=1}^{4}\left(
1+e_{2i-1}\cdot e_{2i}\right) \text{, }s_{2}=\dprod\limits_{i=1}^{4}\left(
1-e_{2i-1}\cdot e_{2i}\right)
\end{equation*}%
are general points of $S_{8}$; their corresponding $8$ by $16$ matrices are:%
\begin{eqnarray*}
P_{0} &=&\left[ I_{8}\left\vert O_{8}\right. \right] \smallskip \\
P_{1} &=&\left[ I_{8}\left\vert J_{4}\right. \right] \smallskip \\
P_{2} &=&\left[ I_{8}\left\vert -J_{4}\right. \right] \text{.}
\end{eqnarray*}%
Let $C$ be the rational normal curve defined by%
\begin{equation*}
C\left( t\right) =\left[ I_{8}\left\vert tJ_{4}\right. \right] \text{.}
\end{equation*}%
We have that $C$ is embedded in $\mathbb{P}^{4}\left( 
%TCIMACRO{\U{2102} }%
%BeginExpansion
\mathbb{C}
%EndExpansion
\right) $, it's contained in $S_{8}$ and%
\begin{equation*}
C\left( 0\right) =P_{0}\text{, }C\left( 1\right) =P_{1}\text{, }C\left(
-1\right) =P_{2}\text{.}
\end{equation*}%
Since%
\begin{equation*}
3\dim _{%
%TCIMACRO{\U{2102} }%
%BeginExpansion
\mathbb{C}
%EndExpansion
}S_{8}+2=86<2^{8-1}-1=127\text{,}
\end{equation*}%
we may apply corollary 3.3 and we get that $\sigma _{3}\left( S_{8}\right) $
hasn't the expected dimension, as desired. \hfill $\square $\medskip

\noindent \textbf{Remark 5.1 }Same argument says that, for all $h=2m$, there
exists a rational normal curve in $S_{h}$ through three points of degree $m$%
.\medskip

\noindent Theorem 5.3 implies that four projectivised tangents spaces to $%
S_{8}$ are always linearly dependent. Hence the following holds:\medskip

\noindent \textbf{Corollary 5.4 }\textit{The variety} $S_{8}$ \textit{is} $4$%
\textit{-defective and }$\delta _{4}=4$.\bigskip

\noindent In the case of $h=7$ we can't apply corollary 5.2. Nevertheless we
have the following:\medskip

\noindent \textbf{Theorem 5.5} \textit{The variety} $S_{7}$ \textit{is} $3$%
\textit{-defective and }$\delta _{3}=5$.\medskip

\noindent \textbf{Proof}. Let $X_{1},X_{2},X_{3}\in S_{7}$ represented in
blocks matrix form and let%
\begin{equation*}
f:SO(14,Q)\rightarrow S_{7}\times S_{7}\times S_{7}
\end{equation*}%
be the function defined by%
\begin{equation*}
f\left( g\right) =\left( g\left( X_{1}\right) ,g\left( X_{2}\right) ,g\left(
X_{3}\right) \right) \text{, for all }g\in SO(14,Q)\text{.}
\end{equation*}%
Taking $g=I_{14}$, the tangent map of $f$ at the point $g$ is:%
\begin{equation*}
df_{I_{14}}:so\left( 14,Q\right) \rightarrow T_{\left(
X_{1},X_{2},X_{3}\right) }\left[ S_{7}\times S_{7}\times S_{7}\right] \text{.%
}
\end{equation*}%
To complete the proof it suffices to find $X_{1}=\left[ I_{7}\left\vert
U_{1}\right. \right] ,$ $X_{2}=\left[ I_{7}\left\vert U_{2}\right. \right] ,$
$X_{3}=\left[ I_{7}\left\vert U_{3}\right. \right] \in S_{7}$ such
that:\smallskip

\noindent 1. the orbit of $\left( X_{1},X_{2},X_{3}\right) $ is open in $%
S_{7}\times S_{7}\times S_{7}$;\smallskip

\noindent 2. $\dim _{%
%TCIMACRO{\U{2102} }%
%BeginExpansion
\mathbb{C}
%EndExpansion
}\left\langle T_{X_{1}}S_{7},T_{X_{2}}S_{7},T_{X_{3}}S_{7}\right\rangle =59$
(we recall that $59$ is the value we got by applying our probabilistic
algorithm at the stage $\left( h,k\right) =\left( 7,3\right) $).\smallskip

\noindent In order that $X_{1},X_{2},X_{3}$ may satisfy the first property,
the rank of the $91$ by $63$ matrix corresponding to $df_{I_{14}}$ has to be
maximum.

\noindent So, we use the first-order Taylor expansion of $f=\left(
f_{1},f_{2},f_{3}\right) $ about $I_{14}$. If%
\begin{equation*}
H=%
\begin{bmatrix}
H_{11} & H_{12} \\ 
H_{21} & H_{22}%
\end{bmatrix}%
\in so\left( 14,Q\right) \text{,}
\end{equation*}%
with $H_{ij}\in M\left( 7,%
%TCIMACRO{\U{2102} }%
%BeginExpansion
\mathbb{C}
%EndExpansion
\right) $, $i,j\in \left\{ 1,2\right\} $, we have that, for $i\in \left\{
1,2,3\right\} $,%
\begin{eqnarray*}
f_{i}\left( I_{14}+H\right) &=&\left[ I_{7}\left\vert \left(
I_{7}+H_{11}^{t}+U_{i}H_{12}^{t}\right) ^{-1}\left( H_{21}^{t}+U_{i}\left(
I_{7}+H_{22}^{t}\right) \right) \right. \right] \smallskip \\
&=&\left[ I_{7}\left\vert
U_{i}+H_{21}^{t}+U_{i}H_{22}^{t}-H_{11}^{t}U_{i}-U_{i}H_{12}^{t}U_{i}+...%
\right. \right] \text{.}
\end{eqnarray*}%
Since $H\in so\left( 14,Q\right) $, it's not hard to show $\left( \cite{Ang}%
\right) $ that, for $i\in \left\{ 1,2,3\right\} $,%
\begin{equation*}
A_{i}=H_{21}^{t}+U_{i}H_{22}^{t}-H_{11}^{t}U_{i}-U_{i}H_{12}^{t}U_{i}^{t}
\end{equation*}%
is a skew symmetric matrix. By computing the jacobian of Pfaffians of size $%
2 $ of $A_{i}$, $i\in \left\{ 0,1,2\right\} $, we get the matrix
corresponding to $df_{I_{14}}$.

\noindent In order to find such points we employed the Macaulay2 software
system, $\cite{Ang}$; in particular $U_{1}=O_{7}$ whereas $U_{2}$ and $U_{3}$
are made of random rational entries. With these choices the above conditions
1. and 2. are satisfied.\hfill $\square $\bigskip

\noindent \textbf{Remark 5.2 }The result of theorem 5.5 agrees with the fact
that the ideal of $\sigma _{2}\left( S_{7}\right) $ is generated in degree $%
4 $, as we can see in $\left( \cite{L-W}\right) $.

\section{Non defective spinor varieties\emph{\protect\medskip }}

\noindent In this section, by using induction, we get our main
result.\medskip

\noindent First of all we have the following:\medskip

\noindent \textbf{Theorem 6.1} \textit{For all }$h\geq 12$, \textit{the
affine tangent spaces to }$S_{h}$ \textit{at}%
\begin{equation*}
P_{0}^{h}=\left[ 
\begin{array}{cc}
I_{12} & O_{12\times (h-12)} \\ 
O_{(h-12)\times 12} & I_{h-12}%
\end{array}%
\left\vert 
\begin{array}{cc}
O_{12} & O_{12\times (h-12)} \\ 
O_{(h-12)\times 12} & O_{h-12}%
\end{array}%
\right. \right] \text{,}
\end{equation*}%
\begin{equation*}
P_{1}^{h}=\left[ 
\begin{array}{cc}
I_{12} & O_{12\times (h-12)} \\ 
O_{(h-12)\times 12} & I_{h-12}%
\end{array}%
\left\vert 
\begin{array}{cc}
J_{6} & O_{12\times (h-12)} \\ 
O_{(h-12)\times 12} & O_{h-12}%
\end{array}%
\right. \right] \text{,}
\end{equation*}%
\begin{equation*}
P_{2}^{h}=\left[ 
\begin{array}{cc}
I_{12} & O_{12\times (h-12)} \\ 
O_{(h-12)\times 12} & I_{h-12}%
\end{array}%
\left\vert 
\begin{array}{cc}
K_{6} & O_{12\times (h-12)} \\ 
O_{(h-12)\times 12} & O_{h-12}%
\end{array}%
\right. \right] \text{,}
\end{equation*}%
\textit{where }$J_{6}$ \textit{is the standard skew symmetric matrix of size}
$12$ \textit{already used before} \textit{and} $K_{6}$ \textit{is the skew
symmetric matrix of size }$12$ \textit{with six diagonal blocks of type}%
\begin{equation*}
\left( 
\begin{array}{cc}
0 & t \\ 
-t & 0%
\end{array}%
\right) \text{, }t\in \left\{ 2,3,..,7\right\} \text{,}
\end{equation*}%
\textit{are linearly independent.}\medskip

\noindent \textbf{Proof}. We proceed by using induction on $h$.\smallskip

\noindent If $h=12$, a slight modification of our probabilistic algorithm in
step $5$ allow us to check the statement.

\noindent Therefore, we assume that the theorem holds for all $h$ such that $%
12\leq h\leq s$, we want to prove it also for $s+1$.

\noindent First of all we remark that $S_{s}$ is embedded in $S_{s+1}$ as
follows:%
\begin{equation}
\left[ I_{s}\left\vert U\right. \right] \in S_{s}\overset{i}{\hookrightarrow 
}\left[ 
\begin{array}{cc}
I_{s} & O_{s\times 1} \\ 
O_{1\times s} & 1%
\end{array}%
\left\vert 
\begin{array}{cc}
U & O_{s\times 1} \\ 
O_{1\times s} & 0%
\end{array}%
\right. \right] =\left[ I_{s+1}\left\vert \widetilde{U}\right. \right] \in
S_{s+1}  \label{2000}
\end{equation}%
where $U\in M\left( s,%
%TCIMACRO{\U{2102} }%
%BeginExpansion
\mathbb{C}
%EndExpansion
\right) $ is\textit{\ }skew symmetric.

\noindent Now, let%
\begin{equation*}
P=\left[ I_{s+1}\left\vert \widetilde{U}\right. \right] =\left[
I_{s+1}\left\vert 
\begin{array}{cccc}
&  &  & y_{1} \\ 
& U &  & \vdots \\ 
&  &  & y_{s} \\ 
-y_{1} & \cdots & -y_{s} & 0%
\end{array}%
\right. \right] \in S_{s+1}\text{,}
\end{equation*}%
with $U=\left\{ u_{ij}\right\} $ skew symmetric of size $s$; we can
parametrize $P$ in $\mathbb{P}^{2^{\left( s+1\right) -1}-1}\left( \mathbb{C}%
\right) $ in such a way that the first coordinates correspond to the
principal sub-Pfaffians of\textit{\ }$U$ and the last one to those of $%
\widetilde{U}$ that involve the last column. Moreover, if $P\in S_{s}$,
then, because of $\left( \ref{2000}\right) $, the affine tangent space to $%
S_{s+1}$ at $i\left( P\right) $ can be represented by the following $\dfrac{%
\left( s+1\right) s}{2}+1\times 2^{\left( s+1\right) -1}$ matrix $M^{s+1}$,
whose blocks form is:%
\begin{equation*}
M^{s+1}=\left( 
\begin{array}{ll}
C_{1} & C_{2}%
\end{array}%
\right)
\end{equation*}%
where%
\begin{equation*}
C_{1}=%
\begin{tabular}{|cccccc|}
\hline
$1$ & $Pf_{2}\left( U\right) $ & $Pf_{4}\left( U\right) $ & $\cdots $ & $%
Pf_{l}\left( U\right) $ & $\cdots $ \\ \hline
$O_{\frac{\left( s-1\right) s}{2}\times 1}$ & $\frac{\partial }{\partial
u_{ij}}Pf_{2}\left( U\right) $ & $\frac{\partial }{\partial u_{ij}}%
Pf_{4}\left( U\right) $ & $\cdots $ & $\frac{\partial }{\partial u_{ij}}%
Pf_{l}\left( U\right) $ & $\cdots $ \\ \hline
&  & $O_{s\times 2^{s-1}}$ &  &  &  \\ \hline
\end{tabular}%
\end{equation*}%
and%
\begin{equation*}
C_{2}=%
\begin{tabular}{|ccc|}
\hline
& $O_{1\times 2^{s-1}}$ &  \\ \hline
& $O_{\frac{\left( s-1\right) s}{2}\times 2^{s-1}}$ &  \\ \hline
$I_{s}$ & \multicolumn{1}{|c}{$A^{s+1}$} & \multicolumn{1}{|c|}{$\ast $} \\ 
\hline
\end{tabular}%
;
\end{equation*}

\noindent $Pf_{l}\left( U\right) $ is the set of the principal sub-Pfaffians
of\textit{\ }$U$ of size $l$, $A^{s+1}$ is the $s\times \dbinom{s}{3}$
matrix made up of \ the derivatives, with respect to $y_{1},...,y_{s}$, of
the principal sub-Pfaffians of\textit{\ }$\widetilde{U}$ of size $4$ that
involve the last column and the entries of $\ast $ are the derivatives, with
respect to $y_{1},...,y_{s}$, of the principal sub-Pfaffians of\textit{\ }$%
\widetilde{U}$ of order $r\geq 6$ that involve the last column.

\noindent We remark that the first two blocks of $C_{1}$%
\begin{equation*}
\begin{tabular}{|cccccc|}
\hline
$1$ & $Pf_{2}\left( U\right) $ & $Pf_{4}\left( U\right) $ & $\cdots $ & $%
Pf_{l}\left( U\right) $ & $\cdots $ \\ \hline
$O_{\frac{\left( s-1\right) s}{2}\times 1}$ & $\frac{\partial }{\partial
u_{ij}}Pf_{2}\left( U\right) $ & $\frac{\partial }{\partial u_{ij}}%
Pf_{4}\left( U\right) $ & $\cdots $ & $\frac{\partial }{\partial u_{ij}}%
Pf_{l}\left( U\right) $ & $\cdots $ \\ \hline
\end{tabular}%
\end{equation*}%
represent the affine tangent space to $S_{s}$ at $P$.

\noindent A direct computation shows that $A^{s+1}$ has the following blocks
structure:%
\begin{equation*}
\left( 
\begin{array}{llll}
D_{1} & D_{2} & \cdots & D_{s-2}%
\end{array}%
\right)
\end{equation*}%
where $D_{i}$'s entries, $i\in \left\{ 1,...,s-2\right\} $, are the
derivatives, with respect to $y_{1},...,y_{s}$, of the principal
sub-Pfaffians of\textit{\ }$\widetilde{U}$ of size $4$ whose first row is
the $i$-th. For our aim, we need only the first four blocks of $A^{s+1}$,
i.e.:%
\begin{equation*}
D_{1}=%
\begin{tabular}{|c|c|c|c|}
\hline
$u_{23}$ $u_{24}$ $\cdots $ $u_{2s}$ & $u_{34}$ $u_{35}$ $\cdots $ $u_{3s}$
& $\cdots $ & $u_{\left( s-1\right) s}$ \\ \hline
$-u_{13}$ $-u_{14}\cdots -u_{1s}$ & $0\cdots \cdots \cdots \cdots 0$ & $%
\cdots $ & $0$ \\ \hline
$u_{12}I_{s-2}$ & $-u_{14}$ $-u_{15}\cdots -u_{1s}$ &  & $0$ \\ 
\cline{2-2}\cline{4-4}
\multicolumn{1}{|l|}{} & $u_{13}I_{s-3}$ & \multicolumn{1}{|l|}{$\ddots $} & 
\multicolumn{1}{|l|}{} \\ 
\multicolumn{1}{|l|}{} & \multicolumn{1}{|l|}{} & \multicolumn{1}{|l|}{$%
\ddots $} & $\vdots $ \\ 
\multicolumn{1}{|l|}{} & \multicolumn{1}{|l|}{} & \multicolumn{1}{|l|}{} & 
\multicolumn{1}{|l|}{} \\ \cline{4-4}
\multicolumn{1}{|l|}{} & \multicolumn{1}{|l|}{} & \multicolumn{1}{|l|}{} & $%
0 $ \\ \cline{4-4}
\multicolumn{1}{|l|}{} & \multicolumn{1}{|l|}{} & \multicolumn{1}{|l|}{} & $%
-u_{1s}$ \\ \cline{4-4}
\multicolumn{1}{|l|}{} & \multicolumn{1}{|l|}{} & \multicolumn{1}{|l|}{} & $%
u_{1\left( s-1\right) }I_{1}$ \\ \hline
\end{tabular}%
\end{equation*}%
\begin{equation*}
D_{2}=%
\begin{tabular}{|c|c|c|c|}
\hline
$0\cdots \cdots \cdots \cdots 0$ & $0\cdots \cdots \cdots \cdots 0$ & $%
\cdots $ & $0$ \\ \hline
$u_{34}$ $u_{35}$ $\cdots $ $u_{3s}$ & $u_{45}$ $u_{46}$ $\cdots $ $u_{4s}$
& $\cdots $ & $u_{\left( s-1\right) s}$ \\ \hline
$-u_{24}$ $-u_{25}\cdots -u_{2s}$ & $0\cdots \cdots \cdots \cdots 0$ &  & $0$
\\ \cline{1-2}\cline{2-2}\cline{4-4}
$u_{23}I_{s-3}$ & $-u_{25}$ $-u_{26}\cdots -u_{2s}$ & $\ddots $ &  \\ 
\cline{2-2}
& $u_{24}I_{s-4}$ & $\ddots $ & $\vdots $ \\ 
&  &  &  \\ \cline{4-4}
&  &  & $0$ \\ \cline{4-4}
&  &  & $-u_{2s}$ \\ \cline{4-4}
&  &  & $u_{2\left( s-1\right) }I_{1}$ \\ \hline
\end{tabular}%
\end{equation*}%
\begin{equation*}
D_{3}=%
\begin{tabular}{|c|c|c|}
\hline
$0\cdots \cdots \cdots \cdots 0$ & $\cdots $ & $0$ \\ \hline
$0\cdots \cdots \cdots \cdots 0$ & $\cdots $ & $0$ \\ \hline
$u_{45}$ $u_{46}$ $\cdots $ $u_{4s}$ & $\cdots $ & $u_{\left( s-1\right) s}$
\\ \hline
$-u_{35}$ $-u_{36}\cdots -u_{3s}$ &  & $0$ \\ \cline{1-1}\cline{3-3}
$u_{34}I_{s-4}$ & $\ddots $ & $\vdots $ \\ 
& $\ddots $ &  \\ \cline{3-3}
&  & $0$ \\ \cline{3-3}
&  & $-u_{3s}$ \\ \cline{3-3}
&  & $u_{3\left( s-1\right) }I_{1}$ \\ \hline
\end{tabular}%
\text{ }D_{4}=%
\begin{tabular}{|c|c|c|}
\hline
$0\cdots \cdots \cdots \cdots 0$ & $\cdots $ & $0$ \\ \hline
$0\cdots \cdots \cdots \cdots 0$ & $\cdots $ & $0$ \\ \hline
$0\cdots \cdots \cdots \cdots 0$ & $\cdots $ & $0$ \\ \hline
$u_{56}$ $u_{57}$ $\cdots $ $u_{5s}$ & $\cdots $ & $u_{\left( s-1\right) s}$
\\ \hline
$-u_{46}$ $-u_{47}\cdots -u_{4s}$ &  & $0$ \\ \cline{1-1}\cline{3-3}
$u_{45}I_{s-5}$ & $\ddots $ & $\vdots $ \\ \cline{3-3}
& $\ddots $ & $0$ \\ \cline{3-3}
&  & $-u_{4s}$ \\ \cline{3-3}
&  & $u_{4\left( s-1\right) }I_{1}$ \\ \hline
\end{tabular}%
\text{.}
\end{equation*}%
So, if instead of a generic skew symmetric $U\in M\left( s,%
%TCIMACRO{\U{2102} }%
%BeginExpansion
\mathbb{C}
%EndExpansion
\right) $, we consider, respectively,%
\begin{equation*}
U_{0}^{s}=\left( 
\begin{array}{cc}
O_{12} & O_{12\times (s-12)} \\ 
O_{(s-12)\times 12} & O_{s-12}%
\end{array}%
\right)
\end{equation*}%
\begin{equation*}
U_{1}^{s}=\left( 
\begin{array}{cc}
J_{6} & O_{12\times (s-12)} \\ 
O_{(s-12)\times 12} & O_{s-12}%
\end{array}%
\right)
\end{equation*}%
\begin{equation*}
U_{2}^{s}=\left( 
\begin{array}{cc}
K_{6} & O_{12\times (s-12)} \\ 
O_{(s-12)\times 12} & O_{s-12}%
\end{array}%
\right)
\end{equation*}%
and we arrange in columns the corresponding $M^{s+1}$ matrices, we get the
span of the affine tangent spaces to $S_{s+1}$ at $P_{0}^{s+1}=i\left(
P_{0}^{s}\right) $, $P_{1}^{s+1}=i\left( P_{1}^{s}\right) $, $%
P_{2}^{s+1}=i\left( P_{2}^{s}\right) $. Reorganizing opportunely the rows,
we can focus our attention on the following $\left[ 3\dfrac{\left(
s+1\right) s}{2}+3\right] \times 2^{\left( s+1\right) -1}$ matrix:%
\begin{equation*}
T^{s+1}=\left( 
\begin{array}{cc}
T^{s} & O_{3\frac{\left( s-1\right) s}{2}+3\times 2^{s-1}} \\ 
O_{3s\times 2^{s-1}} & \Omega%
\end{array}%
\right)
\end{equation*}%
where%
\begin{equation*}
T^{s}=%
\begin{tabular}{|cccccc|}
\hline
$1$ &  & $O_{1\times 2^{s-1}-1}$ &  &  & $\cdots $ \\ \hline
&  & $O_{\frac{\left( s-1\right) s}{2}\times 2^{s-1}-1}$ &  &  & $\cdots $
\\ \hline
$1$ & $Pf_{2}\left( U_{1}^{s}\right) $ & $Pf_{4}\left( U_{1}^{s}\right) $ & $%
\cdots $ & $Pf_{l}\left( U_{1}^{s}\right) $ & $\cdots $ \\ \hline
$O_{\frac{\left( s-1\right) s}{2}\times 1}$ & $\partial Pf_{2\left\vert
U_{1}^{s}\right. }$ & $\partial Pf_{4\left\vert U_{1}^{s}\right. }$ & $%
\cdots $ & $\partial Pf_{l\left\vert U_{1}^{s}\right. }$ & $\cdots $ \\ 
\hline
$1$ & $Pf_{2}\left( U_{2}^{s}\right) $ & $Pf_{4}\left( U_{2}^{s}\right) $ & $%
\cdots $ & $Pf_{l}\left( U_{2}^{s}\right) $ & $\cdots $ \\ \hline
$O_{\frac{\left( s-1\right) s}{2}\times 1}$ & $\partial Pf_{2\left\vert
U_{2}^{s}\right. }$ & $\partial Pf_{4\left\vert U_{2}^{s}\right. }$ & $%
\cdots $ & $\partial Pf_{l\left\vert U_{2}^{s}\right. }$ & $\cdots $ \\ 
\hline
\end{tabular}%
\end{equation*}%
\begin{equation*}
\Omega =%
\begin{tabular}{|l|l|l|}
\hline
$I_{s}$ & $O_{s\times \binom{s}{3}}$ & $\ast _{0}$ \\ \hline
$I_{s}$ & $A_{1}^{s+1}$ & $\ast _{1}$ \\ \hline
$I_{s}$ & $A_{2}^{s+1}$ & $\ast _{2}$ \\ \hline
\end{tabular}%
\text{.}
\end{equation*}%
We want to prove that $T^{s+1}$ has maximum rank, i.e. that%
\begin{equation*}
rankT^{s+1}=3\dfrac{\left( s+1\right) s}{2}+3\text{.}
\end{equation*}%
By induction,%
\begin{equation*}
rankT^{s}=3\dfrac{\left( s-1\right) s}{2}+3\text{,}
\end{equation*}%
being $T^{s}$ the matrix corresponding to the span of the affine tangent
spaces to $S_{s}$ at $P_{0}^{s}$, $P_{1}^{s}$, $P_{2}^{s}$. Then we have
only to prove that 
\begin{equation}
rank\binom{A_{1}^{s+1}}{A_{2}^{s+1}}=2s\text{.}  \label{2500}
\end{equation}

\noindent We remark that $A_{1}^{s+1}=A_{\left\vert U_{1}^{s}\right. }^{s+1}$
and $A_{2}^{s+1}=A_{\left\vert U_{2}^{s}\right. }^{s+1}$; so we consider the
following $2s\times \binom{s}{3}$ blocks matrix:%
\begin{equation*}
\binom{A_{1}^{s+1}}{A_{2}^{s+1}}=\left( 
\begin{array}{llll}
B_{1} & B_{2} & \cdots & B_{s-2}%
\end{array}%
\right)
\end{equation*}%
with $B_{i}=\binom{D_{i\left\vert U_{1}^{s}\right. }}{D_{i\left\vert
U_{2}^{s}\right. }}$, $i\in \left\{ 1,...,s-2\right\} $. In particular we
have that:%
\begin{equation*}
B_{1}=%
\begin{tabular}{|c|c|c|c|c|c|}
\hline
$0\cdots 0$ & $1$ $0\cdots 0$ & $0\cdots 0$ & $1$ $0\cdots 0$ & $\cdots $ & $%
1/0$ \\ \hline
$0\cdots 0$ & $0$ $0\cdots 0$ & $0\cdots 0$ & $0$ $0\cdots 0$ &  & $0$ \\ 
\cline{1-4}\cline{1-4}\cline{6-6}
$I_{s-2}$ & $0$ $0\cdots 0$ & $0\cdots 0$ & $0$ $0\cdots 0$ &  &  \\ 
\cline{2-4}\cline{3-4}
& $O_{s-3}$ & $0\cdots 0$ & $0$ $0\cdots 0$ &  &  \\ \cline{3-4}\cline{4-4}
&  & $O_{s-4}$ & $0$ $0\cdots 0$ &  & $\vdots $ \\ \cline{4-4}
&  &  & $O_{s-5}$ &  &  \\ 
&  &  &  & $\ddots $ &  \\ 
&  &  &  &  &  \\ \cline{6-6}
&  &  &  &  & $O_{1}$ \\ \hline
$0\cdots 0$ & $3$ $0\cdots 0$ & $0\cdots 0$ & $4$ $0\cdots 0$ & $\cdots $ & $%
7/0$ \\ \hline
$0\cdots 0$ & $0$ $0\cdots 0$ & $0\cdots 0$ & $0$ $0\cdots 0$ &  & $0$ \\ 
\cline{1-4}\cline{6-6}
$2I_{s-2}$ & $0$ $0\cdots 0$ & $0\cdots 0$ & $0$ $0\cdots 0$ &  &  \\ 
\cline{2-4}
& $O_{s-3}$ & $0\cdots 0$ & $0$ $0\cdots 0$ &  &  \\ \cline{3-4}
&  & $O_{s-4}$ & $0$ $0\cdots 0$ &  & $\vdots $ \\ \cline{4-4}
&  &  & $O_{s-5}$ &  &  \\ 
&  &  &  & $\ddots $ &  \\ 
&  &  &  &  &  \\ \cline{6-6}
&  &  &  &  & $O_{1}$ \\ \hline
\end{tabular}%
\text{,}
\end{equation*}%
\begin{equation*}
B_{2}=%
\begin{tabular}{|c|c|c|c|c|}
\hline
$0$ $0\cdots 0$ & $0\cdots 0$ & $0$ $0\cdots 0$ & $\cdots $ & $0$ \\ \hline
$1$ $0\cdots 0$ & $0\cdots 0$ & $1$ $0\cdots 0$ &  & $1/0$ \\ 
\cline{1-1}\cline{1-3}\cline{2-3}\cline{5-5}
$0$ $0\cdots 0$ & $0\cdots 0$ & $0$ $0\cdots 0$ &  & $0$ \\ 
\cline{1-3}\cline{2-3}\cline{5-5}
$O_{s-3}$ & $0\cdots 0$ & $0$ $0\cdots 0$ &  &  \\ \cline{2-3}\cline{3-3}
& $O_{s-4}$ & $0$ $0\cdots 0$ &  &  \\ \cline{3-3}
&  & $O_{s-5}$ &  & $\vdots $ \\ 
&  &  & $\ddots $ &  \\ 
&  &  &  &  \\ \cline{5-5}
&  &  &  & $O_{1}$ \\ \hline
$0$ $0\cdots 0$ & $0\cdots 0$ & $0$ $0\cdots 0$ & $\cdots $ & $0$ \\ \hline
$3$ $0\cdots 0$ & $0\cdots 0$ & $4$ $0\cdots 0$ &  & $7/0$ \\ 
\cline{1-3}\cline{5-5}
$0$ $0\cdots 0$ & $0\cdots 0$ & $0$ $0\cdots 0$ &  & $0$ \\ 
\cline{1-3}\cline{5-5}
$O_{s-3}$ & $0\cdots 0$ & $0$ $0\cdots 0$ &  &  \\ \cline{2-3}
& $O_{s-4}$ & $0$ $0\cdots 0$ &  &  \\ \cline{3-3}
&  & $O_{s-5}$ &  & $\vdots $ \\ 
&  &  & $\ddots $ &  \\ 
&  &  &  &  \\ \cline{5-5}
&  &  &  & $O_{1}$ \\ \hline
\end{tabular}%
\text{,}
\end{equation*}%
\begin{equation*}
B_{3}=%
\begin{tabular}{|c|c|c|c|}
\hline
$0\cdots 0$ & $0$ $0\cdots 0$ & $\cdots $ & $0$ \\ \hline
$0\cdots 0$ & $0$ $0\cdots 0$ &  & $0$ \\ 
\cline{1-1}\cline{1-2}\cline{2-2}\cline{4-4}
$0\cdots 0$ & $1$ $0\cdots 0$ &  & $1/0$ \\ 
\cline{1-1}\cline{1-2}\cline{2-2}\cline{4-4}
$0\cdots 0$ & $0$ $0\cdots 0$ &  & $0$ \\ \cline{1-2}\cline{2-2}\cline{4-4}
$I_{s-4}$ & $0$ $0\cdots 0$ &  &  \\ \cline{2-2}
& $O_{s-5}$ &  & $\vdots $ \\ 
&  & $\ddots $ &  \\ 
&  &  &  \\ \cline{4-4}
&  &  & $O_{1}$ \\ \hline
$0\cdots 0$ & $0$ $0\cdots 0$ & $\cdots $ & $0$ \\ \hline
$0\cdots 0$ & $0$ $0\cdots 0$ &  & $0$ \\ \cline{1-2}\cline{4-4}
$0\cdots 0$ & $4$ $0\cdots 0$ &  & $7/0$ \\ \cline{1-2}\cline{4-4}
$0\cdots 0$ & $0$ $0\cdots 0$ &  & $0$ \\ \cline{1-2}\cline{4-4}
$3I_{s-4}$ & $0$ $0\cdots 0$ &  &  \\ \cline{2-2}
& $O_{s-5}$ &  & $\vdots $ \\ 
&  & $\ddots $ &  \\ 
&  &  &  \\ \cline{4-4}
&  &  & $O_{1}$ \\ \hline
\end{tabular}%
\text{, }B_{4}=%
\begin{tabular}{|c|c|c|}
\hline
$0$ $0\cdots 0$ & $\cdots $ & $0$ \\ \hline
$0$ $0\cdots 0$ &  & $0$ \\ \cline{1-1}\cline{1-1}\cline{3-3}
$0$ $0\cdots 0$ &  & $0$ \\ \cline{1-1}\cline{1-1}\cline{3-3}
$1$ $0\cdots 0$ &  & $1/0$ \\ \cline{1-1}\cline{1-1}\cline{3-3}
$0$ $0\cdots 0$ &  & $0$ \\ \cline{1-1}\cline{3-3}
$O_{s-5}$ &  &  \\ 
& $\ddots $ & $\vdots $ \\ 
&  &  \\ \cline{3-3}
&  & $O_{1}$ \\ \hline
$0$ $0\cdots 0$ & $\cdots $ & $0$ \\ \hline
$0$ $0\cdots 0$ &  & $0$ \\ \cline{1-1}\cline{3-3}
$0$ $0\cdots 0$ &  & $0$ \\ \cline{1-1}\cline{3-3}
$4$ $0\cdots 0$ &  & $7/0$ \\ \cline{1-1}\cline{3-3}
$0$ $0\cdots 0$ &  & $0$ \\ \cline{1-1}\cline{3-3}
$O_{s-5}$ &  &  \\ 
& $\ddots $ & $\vdots $ \\ 
&  &  \\ \cline{3-3}
&  & $O_{1}$ \\ \hline
\end{tabular}%
\text{.}
\end{equation*}%
We observe that in the case of $s=12$ we consider the element before $/$,
otherwise the element after.

\noindent By the Gauss elimination algorithm, the blocks $B_{1}$, $B_{2}$, $%
B_{3}$ and $B_{4}$ become, respectively:%
\begin{equation*}
\overline{B_{1}}=%
\begin{tabular}{|c|c|c|c|c|c|}
\hline
& $0$ $0\cdots 0$ & $0\cdots 0$ & $0$ $0\cdots 0$ & $\cdots $ & $0$ \\ 
\cline{2-6}\cline{3-3}\cline{6-6}
$I_{s-2}$ &  & $0\cdots 0$ & $0$ $0\cdots 0$ &  &  \\ \cline{3-4}
& $O_{s-3}$ &  & $0$ $0\cdots 0$ &  &  \\ \cline{4-4}
&  & $O_{s-4}$ &  &  &  \\ 
&  &  & $O_{s-5}$ &  & $\vdots $ \\ 
&  &  &  & $\ddots $ &  \\ \cline{6-6}
&  &  &  &  & $O_{1}$ \\ \hline
$0\cdots 0$ & $1$ $0\cdots 0$ & $0\cdots 0$ & $1$ $0\cdots 0$ &  & $1/0$ \\ 
\hline
$0\cdots 0$ & $0$ $0\cdots 0$ & $0\cdots 0$ & $0$ $0\cdots 0$ &  & $0$ \\ 
\hline
$0\cdots 0$ & $0$ $0\cdots 0$ & $0\cdots 0$ & $1$ $0\cdots 0$ & $\cdots $ & $%
4/0$ \\ \hline
$0\cdots 0$ & $0$ $0\cdots 0$ & $0\cdots 0$ & $0$ $0\cdots 0$ &  & $0$ \\ 
\cline{1-4}\cline{6-6}
& $0$ $0\cdots 0$ & $0\cdots 0$ & $0$ $0\cdots 0$ &  &  \\ \cline{2-4}
$O_{s-2}$ &  & $0\cdots 0$ & $0$ $0\cdots 0$ &  &  \\ \cline{3-4}
& $O_{s-3}$ &  & $0$ $0\cdots 0$ &  & $\vdots $ \\ \cline{4-4}
&  & $O_{s-4}$ &  &  &  \\ 
&  &  & $O_{s-5}$ &  &  \\ \cline{6-6}
&  &  &  & $\ddots $ & $0$ \\ \cline{6-6}
&  &  &  &  & $O_{1}$ \\ \hline
\end{tabular}%
\end{equation*}%
\begin{equation*}
\overline{B_{2}}=%
\begin{tabular}{|c|c|c|c|c|}
\hline
$0$ $0\cdots 0$ & $0\cdots 0$ & $0$ $0\cdots 0$ & $\cdots $ & $0$ \\ \hline
& $0\cdots 0$ & $0$ $0\cdots 0$ &  &  \\ \cline{2-3}\cline{2-3}
$O_{s-3}$ &  & $0$ $0\cdots 0$ &  &  \\ \cline{3-3}\cline{3-3}
& $O_{s-4}$ &  &  & $\vdots $ \\ 
&  & $O_{s-5}$ &  &  \\ 
&  &  & $\ddots $ &  \\ \cline{5-5}
&  &  &  & $O_{1}$ \\ \hline
$0$ $0\cdots 0$ & $0\cdots 0$ & $0$ $0\cdots 0$ &  & $0$ \\ \hline
$1$ $0\cdots 0$ & $0\cdots 0$ & $1$ $0\cdots 0$ & $\cdots $ & $1/0$ \\ \hline
$0$ $0\cdots 0$ & $0\cdots 0$ & $0$ $0\cdots 0$ & $\cdots $ & $0$ \\ \hline
$0$ $0\cdots 0$ & $0\cdots 0$ & $1$ $0\cdots 0$ &  & $4/0$ \\ 
\cline{1-3}\cline{5-5}
$0$ $0\cdots 0$ & $0\cdots 0$ & $0$ $0\cdots 0$ &  & $0$ \\ 
\cline{1-3}\cline{5-5}
& $0\cdots 0$ & $0$ $0\cdots 0$ &  &  \\ \cline{2-3}
$O_{s-3}$ &  & $0$ $0\cdots 0$ &  &  \\ \cline{3-3}
& $O_{s-4}$ &  &  & $\vdots $ \\ 
&  & $O_{s-5}$ &  &  \\ \cline{5-5}
&  &  & $\ddots $ & $0$ \\ \cline{5-5}
&  &  &  & $O_{1}$ \\ \hline
\end{tabular}%
\end{equation*}%
\begin{equation*}
\overline{B_{3}}=%
\begin{tabular}{|c|c|c|c|}
\hline
$0\cdots 0$ & $1$ $0\cdots 0$ & $\cdots $ & $1/0$ \\ \hline
$0\cdots 0$ & $0$ $0\cdots 0$ &  & $0$ \\ 
\cline{1-1}\cline{1-2}\cline{2-2}\cline{4-4}
& $0$ $0\cdots 0$ &  &  \\ \cline{2-2}\cline{2-2}
$I_{s-4}$ &  &  & $\vdots $ \\ 
& $O_{s-5}$ &  &  \\ \cline{4-4}
&  & $\ddots $ & $0$ \\ \cline{4-4}
&  &  & $O_{1}$ \\ \cline{1-2}\cline{4-4}
$0\cdots 0$ & $0$ $0\cdots 0$ &  & $0$ \\ \cline{1-2}\cline{4-4}
$0\cdots 0$ & $0$ $0\cdots 0$ &  & $0$ \\ \hline
$0\cdots 0$ & $0$ $0\cdots 0$ & $\cdots $ & $0$ \\ \hline
$0\cdots 0$ & $0$ $0\cdots 0$ &  & $0$ \\ \cline{1-2}\cline{4-4}
$0\cdots 0$ & $2$ $0\cdots 0$ &  & $5/0$ \\ \cline{1-2}\cline{4-4}
$0\cdots 0$ & $0$ $0\cdots 0$ &  & $0$ \\ \cline{1-2}\cline{4-4}
& $0$ $0\cdots 0$ &  &  \\ \cline{2-2}
$I_{s-4}$ &  &  & $\vdots $ \\ 
& $O_{s-5}$ &  &  \\ \cline{4-4}
&  & $\ddots $ & $0$ \\ \cline{4-4}
&  &  & $O_{1}$ \\ \hline
\end{tabular}%
\text{, }\overline{B_{4}}=%
\begin{tabular}{|c|c|c|}
\hline
$0$ $0\cdots 0$ & $\cdots $ & $0$ \\ \hline
$1$ $0\cdots 0$ &  & $1/0$ \\ \cline{1-1}\cline{1-1}\cline{3-3}
$0$ $0\cdots 0$ &  & $0$ \\ \cline{1-1}\cline{1-1}\cline{3-3}
&  & $\vdots $ \\ 
$O_{s-5}$ &  &  \\ \cline{3-3}
& $\ddots $ & $0$ \\ \cline{3-3}
&  & $O_{1}$ \\ \cline{1-1}\cline{3-3}
$0$ $0\cdots 0$ &  & $0$ \\ \cline{1-1}\cline{3-3}
$0$ $0\cdots 0$ &  & $0$ \\ \hline
$0$ $0\cdots 0$ & $\cdots $ & $0$ \\ \hline
$0$ $0\cdots 0$ &  & $0$ \\ \cline{1-1}\cline{3-3}
$0$ $0\cdots 0$ &  & $0$ \\ \cline{1-1}\cline{3-3}
$2$ $0\cdots 0$ &  & $5/0$ \\ \cline{1-1}\cline{3-3}
$0$ $0\cdots 0$ &  & $0$ \\ \cline{1-1}\cline{3-3}
&  &  \\ 
$O_{s-5}$ &  & $\vdots $ \\ \cline{3-3}
& $\ddots $ & $0$ \\ \cline{3-3}
&  & $O_{1}$ \\ \hline
\end{tabular}%
\text{.}
\end{equation*}%
Now it's easy to check that $\left( \ref{2500}\right) $ holds, as desired.
\hfill $\square $\medskip

\noindent As a consequence we get immediately:\medskip

\noindent \textbf{Theorem 6.2} \textit{For all }$h\geq 12$, $\sigma
_{3}\left( S_{h}\right) $ \textit{has the expected dimension.}

\noindent {\small Author's address:}

\noindent Dipartimento di Matematica "Ulisse Dini"

\noindent Viale G. B. Morgagni, 67/a - 50134 Firenze, Italy

\noindent email: elena.angelini@math.unifi.it

\end{document}